\magnification=1200
\baselineskip=14pt
%
%
\font\gen=cmr10 at 10pt
\font\genb=cmbx10 at 10 pt
\font\genit=cmti10
\font\genhead=cmbx10 at 12 pt
\font\gensmall=cmr10 at 9pt
\font\gensmallB=cmbx10 at 9pt
\gen
%
%
\def\qed{${\vcenter{\vbox{\hrule height .4pt
           \hbox{\vrule width .4pt height 4pt
            \kern 4pt \vrule width .4pt}
             \hrule height .4pt}}}$}
\def\mqed{{\vcenter{\vbox{\hrule height .4pt
           \hbox{\vrule width .4pt height 4pt
            \kern 4pt \vrule width .4pt}
             \hrule height .4pt}}}}
\font\bbb=msbm10 at 10pt
\def\NN{\hbox{\bbb N}}
\def\RR{\hbox{\bbb R}}

\def\ZZ{\hbox{\bbb Z}}

\def\QQ{\hbox{\bbb Q}}

\newcount\refCount
\def\newref#1 {\advance\refCount by 1
\expandafter\edef\csname#1\endcsname{\the\refCount}}
\newref BOMB  
\newref BOPI  
\newref GAVO  
\newref GWKP  %
\newref HEAT  
\newref KHOV  
\newref LANG  
\newref PILA  
\newref PILC  
\newref PILD  
\newref SHID  
%
%
\centerline{\genhead 
Note on the rational points of a pfaff curve\/}
\medskip
\bigskip
\noindent
\centerline{\gen Jonathan Pila}
\bigskip
\medskip
{
\baselineskip=11 pt
\midinsert
\centerline{\gensmallB Abstract\/} 
\narrower
\narrower
\smallskip
\noindent
{\gensmall 
Let  $X\subset\RR^2$ be the graph of a pfaffian function $f$ in the
sense  of Khovanskii. Suppose that $X$ is nonalgebraic. This note
gives an estimate for the number of rational points on $X$ of height
$\le H$; the estimate is uniform in the order and degree of $f$.}
\bigskip
\noindent
{{\gensmallB 2000 Mathematics Subject Classification:\/}
{\gensmall 11D99 (11J99)\/}}

\endinsert
}
\bigskip
\centerline{\genb 1. Introduction\/} 
\medskip

In [\PILA] and [\PILD] I have studied the distribution of rational
points on the graph $X$ of a transcendental real analytic
function $f$ on a {\genit compact\/} interval. I showed that the
number of rational points of $X$ of {\genit height\/} (see 1.2 below)
$\le H$ was $O_{f,\epsilon}(H^\epsilon)$ for all positive $\epsilon$. 

\smallskip

Suppose that $X$ is the graph of a function analytic
on a {\genit noncompact} domain, such as $\RR$ or $\RR^+$.
To bound the number of rational points of height $\le
H$ on $X$ requires controlling the implied constant in the above
estimate over the enlarging intervals $[-H,H]$.

\smallskip

In the estimate in [\PILD], the implied constant depends on a bound
for the number of solutions of an algebraic equation in $P(x,f)$, where
$P\in\RR[x,y]$ is a polynomial (of degree depending on $\epsilon$), as
well as a bound for the number of zeros of derivatives of $f$ (of
order depending on $\epsilon$). In general these quantities may not
behave at all well over different intervals. 

\smallskip

However these numbers are globally bounded for the
socalled {\genit pfaffian functions\/} ([\KHOV, \GAVO]; see 1.1 below),
indeed bounded uniformly in terms of the {\genit order\/} and {\genit
degree\/} of the function (see 1.1).  For this class of functions,
a uniform estimate on the number of rational points of bounded
height may be obtained by adapting the methods of [\BOPI, \PILA, \PILD].

\medskip
\noindent
{\genb Definition 1.1.\/} ([\GAVO, 2.1]) Let $U\subset\RR^n$ be an open
domain. A {\genit pfaffian chain\/} of order $r\ge 0$ and degree
$\alpha\ge 1$ in $U$ is a sequence of real analytic functions
$f_1,\ldots,f_r$ in $U$ satisfying differential equations
$$
df_j=\sum_{i=1}^n g_{ij}({\bf x}, f_1({\bf x}), f_2({\bf x}),\ldots,
f_j({\bf x}))dx_i
$$ 
for $j=1,\ldots,r$,
where ${\bf x}=(x_1,\ldots,x_n)$ and $g_{ij}\in\RR[x_1,\ldots,x_n,
y_1,\ldots,y_r]$ of degree $\le \alpha$. A function
$f$ on $U$ is called a {\genit pfaffian function} of order $r$ and
degree
$(\alpha,\beta)$ if
$f({\bf
x})=P({\bf x}, f_1({\bf x}),\ldots,f_r({\bf x}))$, where $P$ is a
polynomial of degree at most $\beta\ge 1$. 

\bigbreak

The usual elementary functions $e^x, \log x$ (but not $\sin x$ on all
$\RR$), algebraic functions, combinations and compositions of these
are pfaffian functions: see [\KHOV, \GAVO]. In this paper always $n=1$,
so ${\bf x}=x$. A {\genit pfaff curve\/} $X$ is the graph of a pfaffian
function
$f$ on some connected subset of its domain. The order and degree of $X$
will be taken to be the order and degree of $f$.

\medskip
\noindent
{\genb Definition 1.2.\/} For a point $P=(a_1/b_1,
a_2/b_2,\ldots,a_n/b_n)\in\QQ^n$, where $a_j,b_j\in\ZZ, b_j\ge 1$ and
$(a_j,b_j)=1$ for all $j=1,2,\ldots,n$, define the {\genit height\/}
$H(P)= \max\{|a_j|, b_j\}$.
Note this is not the projective height. If
$X\subset\RR^n$ let $X(\QQ)=X\cap\QQ^n$ and $X(\QQ,H)$ the subset of
points $P$ with $H(P)\le H$. Finally put
$$
N(X,H)=\#X(\QQ,H)=\#\{P\in X(\QQ), H(P)\le H\}.
$$

\medskip
\noindent
{\genb Theorem 1.3.\/} {\genit
There is an explicit function $c(r,\alpha,\beta)$ with the following
property. Suppose $X$ is a nonalgebraic pfaff curve of order $r$ and
degree $(\alpha,\beta)$. Let $H\ge c(r,\alpha,\beta)$. Then
$$
N(X,H)\le \exp\big(5\sqrt{\log H}\big).
$$

}
\medskip

Now in certain cases where the polynomials defining the
chain have rational (or
algebraic) coefficients, results in transcendence theory show that the
number of {\genit algebraic\/} points of $X$ is {\genit finite}, indeed
explicitly bounded (see e.g. [\LANG, \SHID]). On the other hand, the
example $X=\{(x,y): y=2^x, x\in\RR\}$ shows that the set $X(\QQ)$ is not
finite in general. For many $X$, e.g. the graph of
$y=e^{e^x}$, finiteness is unknown. 

\smallskip

For the example $X=\{(x,y): y=2^x, x\in\RR\}$ of course $N(X,H)=O(\log
H)$. I do not know examples where the growth of $N(X,H)$
is faster than this, so the above bound might be very far from the
truth. Note however that elementary considerations do not suffice to
establish better bounds on $N(X,H)$ for, e.g.,
$X=\{(x,y), y=\log\log(e^{e^x}+e^x)\}$,
for which finiteness of $X(\QQ)$ is presumably expected.

\smallskip

The methods herein are also applicable to algebraic curves: indeed the
fact that pfaffian functions have finiteness properties analogous to
algebraic functions was the impetus for applying those
methods to them. Since algebraic functions are pfaffian ([\KHOV]), it
is appropriate to record here the following improvement
to the result obtained in [\PILD]. 

\medskip
\noindent
{\genb Theorem 1.4.\/} {\genit Let $b,c\ge 2$ be integers and $H\ge
3$. Let $F(x,y)\in\RR[x,y]$ be irreducible of bidegree $(b,c)$. 
Let $d=\max(b,c)$ and $X=\{(x,y)\in\RR^2, F(x,y)=0\}$. Then
$$
N(X,H)\le (6d)^{10}4^d\,H^{2/d}\,\big(\log H\big)^5.
$$

}

\medskip

The improvement over the result of [\PILD] is that the exponent
of $\log H$ is here independent of $d$; 
the exponent $2/d$ of $H$ is best possible.
I refer to [\PILD] for
discussion of related results ([\HEAT, \BOMB]).

\bigbreak
\centerline{\genb 2. The Main Lemma\/}
\medskip
\noindent

Let $M=\{x^hy^k: (h,k)\in J\}$ be a finite set of monomials in the
indeterminates $x,y$. Put
$$
D=\#M,\,\,\,\,\,\,\, 
R=\sum_{(h,k)\in J}(h+k),\,\,\,\,\,\,\, 
s=\max_{(h,k)\in J}(h),\,\,\,\,\,\,\, 
t=\max_{(h,k)\in J}(k),\,\,\,\,\,\,\,
S=D\big(s+t\big),
$$
$$
\rho={2R\over D(D-1)},\,\,\,\,\,\,\,
\sigma={2S\over D(D-1)},\,\,\,\,\,\,\,
C=\big(D!D^R\big)^{2/D(D-1)}+1.
$$
Note that $S\ge R$ for any $M$.
If $Y$ is a plane algebraic curve defined by $G(x,y)=0$, 
say $Y$ is {\genit defined in $M$\/} if all the monomials appearing in
$G$ belong to $M$.

\medskip
\noindent
{\genb Lemma 2.1.\/} {\genit
Let $M$ be a set of monomials with $D\ge 2$ and $S\ge 2R$.
Let $H\ge 1, L\ge 1/H^2$ and $I$ a closed interval of length $\le L$.
Let $f\in C^D(I)$ with $|f'|\le 1$ and $f^{(j)}$ either nonvanishing in
the interior of $I$ or identically zero for $j=1,2,\ldots, D$.
Let $X$ be the graph of $y=f(x)$ on $I$.
Then $X(\QQ, H)$ is contained in the union of at most
$$
\big(4\,C\,D\,4^{1/\rho}+2\big)\,L^\rho\,H^\sigma
$$
real algebraic curves defined in $M$.
}
\medskip
\noindent
{\genb Proof.\/} 
Fix $M, H$. If $f$ is a function satisfying the hypotheses on some
interval $I$, and $X$ is the graph of $f$ on $I$, then the set
$X(\QQ,H)$ is contained in some minimal number $G(f,I)$ of
algeraic curves of degree $\le d$;
Let the $G(L)$ be the maximum of $G(f,I)$ over all intervals and
functions satisfying the hypotheses.

\smallskip

Now suppose $f$ is such a function on an interval $I=[a,b]$, and $A\ge
1$. An equation $f^{(2)}(x)=\pm 2A\,L^{-1}$ has at most one
solution in the interior $I$, unless it is satisfied identically.
Suppose $c$ is a solution. Since
$f^{(2)}, f^{(3)}$ are onesigned throughout $I$, it follows that
$|f^{(2)}(x)|\le  2\,A^{2/(D-1)}\,L^{-1}$ in either $[a,c]$ or $[c,b]$,
and $|f^{(2)}(x)|\ge  2\,A^{2/(D-1)}\,L^{-1}$ in (respectively) 
either $[b,c]$ or $[a,c]$. Now an interval with the latter condition
has length 
$\le 2\,A^{1/(D-1)}$ by [\PILD, 2.6] (or [\BOPI, Lemma 7]) applied with
$A=A^{2/(D-1)}$.

\smallskip

Continuing to split the interval at points where
$f^{(\kappa)}=\kappa!\,A^{\kappa/(D-1)}\,L^{1-\kappa},
\kappa=2,\ldots,D$ yields a (possibly empty) subinterval $[s,t]$ in
which 
$|f^{(\kappa)}|\le\kappa!\,A^{\kappa/(D-1)}\,L^{1-\kappa}$ for all
$\kappa=1,2,3,\ldots,D$, while the intervals $[a,s], [t,b]$ comprise
$\le D$ subintervals of length $\le 2A^{-1/(D-1)}L$
(by [\PILD, 2.6] (or [\BOPI, Lemma 7]) applied with
$A=A^{\kappa/(D-1)}$), and so have length at most $2DA^{-1/(D-1)}L$. 
(If $[s,t]$ is empty take $s=t=b$.)

\smallskip

On $[s,t]$, the points of height $\le H$ lie on at most
$C\,A^{1/(D-1)}\,H^\sigma\,L^\rho$
curves in $M$ by [\PILD, 2.4].
Therefore the function $G(L)$ satisfies the recurrence
$$
G(L)\le C\,A^{1/(D-1)}\,H^{\sigma}L^\rho+2\,G(\lambda\, L)
$$
when $L\ge 1/H^2$, where $\lambda=2\,D\,A^{-1/(D-1)}$.
Thus, provided $\lambda^{n-1}L\ge 1/H^2$,
$$
G(L)\le
C\,A^{1/(D-1)}\,H^\sigma\,L^\rho\,\big(1+2\lambda^\rho+\ldots+
(2\lambda^\rho)^{n-1}\big)\,+\, 2^n\,G(\lambda^n L).
$$
Choose $A$ such that $2\lambda^\rho=1/2$, that is
$A^{1/(D-1)}=2\,D\,4^{1/\rho}$
(so $A\ge 1$) and choose $n$ such that
$$
{\lambda\over LH^2}\le\lambda^n<{1\over LH^2}.
$$
Then $G(\lambda^n L)\le 1$, while 
$$
2^n=\lambda^{-n\rho/2}\le \Big({LH^2\over \lambda}\Big)^{\rho/2}
=2\,\big(L\,H^2\big)^{\rho/2}\le 2L^\rho H^\sigma.
$$
Therefore $G(L)\le \big(4CD\, 4^{1/\alpha}+2\big)\,H^\rho\,L^\sigma$
as required.\ \qed
\bigskip
\bigbreak
\centerline{\genb 3. Pfaff curves\/}
\medskip
\noindent

Since a pfaffian function of order $r=0$ is a polynomial, to which
Theorem 1.3 is inapplicable, it is convenient now to assume $r\ge 1$.

\medskip
\noindent
{\genb Proposition 3.1.\/} {\genit
Let $f_1,\ldots,f_r$ be a pfaffian chain of order $r\ge 1$ and degree
$\alpha$ on an open domain $U\subset\RR$, and $f$ a pfaffian function
on $U$ having this chain and degree $(\alpha,\beta)$.

\smallskip

(a). Let $k\in\NN$. Then $f^{(k)}$ is a pfaffian function with the same
chain as $f$ (so of order $r$) and degree 
$\big(\alpha, \beta+k(\alpha-1)\big)$.

\smallskip

(b). Let $P(x,y)$ be a polynomial of degree $d$. Suppose $f$
is not algebraic. Then the equation
$P(x,f(x))=0$ has at most 
$$
2^{r(r-1)/2}\,d\beta\,\big(r\alpha+d\beta\big)^r
$$ 
solutions.

\smallskip

(c). Let $V\subset U$ be an open set on which $f'\neq 0$ and $k\ge 1$.
Then  on $V$ there is an inverse function $g$ of $f$, and the number of
zeros of $g^{(k)}$ on $V$ is at most
$$
2^{r(r-1)/2}\,\big((k-1)(\beta+k(\alpha-1))\big)\,
\Big(r\alpha+\big((k-1)(\beta+k(\alpha-1))\big)\Big)^r.
$$

}

\medskip
\noindent
{\genb Proof.\/} 
Part (a) is by [\GAVO, 2.5].
For part (b), observe that $P(x,f(x))$ is a pfaffian function of order
$r\ge 1$ and degree $(\alpha,d\beta)$. Since $f$ is not algebraic, all
the solutions are nondegenerate and the result is in [\GAVO, 3.1].

\smallskip

Part (c). By differentiating the relation $g(f(x))=x$ and simple
induction, for $k\ge 1$,
$$
g^{(k)}(y)={Q_k(f^{(1)}, f^{(2)},\ldots, f^{(k)})\over (f'(x))^{2k-1}}
$$
where $Q_k(z_1, z_2,\ldots,z_k)$ is a polynomial of degree
$\gamma_k=k-1$. Since
$f^{(j)}$ are pfaffian functions with the same chain, the function
$Q_k(f^{(1)}, f^{(2)},\ldots, f^{(k)})$ is a pfaffian function of order
$r$ and degree $(\alpha,
\gamma_k(\beta+k(\alpha-1)))$. The statement now follows from (b).\ \qed

\medskip
\noindent
{\genb Proof of 1.3.\/} Let $d\ge 2$ and let $M=M(d)$ be the set
of monomials of degree $d$ in $x,y$. Then, [\PILD] elementarily,
$$
D={d(d-1)\over 2},\,\,\,\,\,\,\, \rho={8\over 3(d+3)},\,\,\,\,\,\,\,
\sigma=3\rho,\,\,\,\,\,\,\, C\le 6.
$$

Subdivide the connected domain $U$ into at most
$$
2.2^{r(r-1)/2}(\beta+\alpha-1)(r\alpha+\beta+\alpha-1)+1 \le 
2^{1+r(r-1)/2}\big((r+1)(\alpha+\beta)\big)^{r+1}
$$
intervals on which $f'\le -1$, $-1\le f'\le 1$ or $f'\ge 1$, and then
divide further into subintervals on which the inverse $g$ has
nonvanishing derivatives up to order $D$ in the first and third case,
or $f$ has nonvanishing derivatives up to order $D$ in the second case.
The total number of intervals is at most
$$
2^{1+r(r-1)/2}\big((r+1)(\alpha+\beta)\big)^{r+1}\,
D^2\,2^{r(r-1)}\,(\beta+D(\alpha-1))(r+D(\beta+D(\alpha-1)))^r.
$$ 

Intersecting with the interval $[-H,H]$ of the appropriate
axis, the intervals of length $\le 2H$. By 2.1, in each interval the
points of $X(\QQ,H)$ lie on at most
$$
(24D4^{1/\rho}+2) (2H)^\rho H^{3\rho}\le 6d^2 4^{1/\rho} 2^\rho
$$
real algebraic curves of degree $d$; the number of points of $X$
on a curve of degree $d$ is at most
$$
2^{r(r-1)/2}\,d\beta\,(r\alpha+d\beta)^r.
$$
Combining these estimates yields
$$
N(X,H)\le c'(r,\alpha,\beta,d,D)\, 4^{3(d+3)/8} H^{32/(3(d+3)}
$$
where $t=3(d+3)/8$. Choose $d$ so that
$t$ is as near as possible to (and so within $1/2$ of) $\sqrt{4\log
H/\log 4}$. Then $4\sqrt{\log 4}<5$ and noting that
$d,D$ appear polynomially in $c'$ completes the proof.\ \qed

\bigbreak

\noindent
{\genb Remarks 3.2.\/}

\smallskip

1. Note that the constant `5' in 1.3 can be improved by further
optimizing the proof. However, a bound of the shape $\exp(c\sqrt{\log
H})$ seems to be the best obtainable by the present method.

\smallskip

2. A result can be formulated for any real analytic (or even smooth)
function $f$ with suitable finiteness properties (zeros of derivatives,
derivatives of the inverse, and algebraic relations). An example of
such a function that is not pfaffian is exhibited in [\GWKP]. (Indeed
the given example $e^x+\sin x$ does not belong to any {\genit
$o$-minimal structure\/}: see [\GWKP]). 

\smallskip

3. I expect a similar result should hold in higher dimensions for pfaff
manifolds: that is, a uniform (in `complexity') $H^\epsilon$ bound for
rational points that do not lie on some semialgebraic subset of positive
dimension (cf the conjectures for subanalytic sets made in
[\PILC, \PILD]).  A similar result should hold for sets definable in an
$o$-minimal structure. 

\bigbreak
\centerline{\genb 4. Algebraic curves}
\medskip
\noindent

For integers $\beta,\gamma\ge 2$ let 
$$
M(\beta,\gamma)=\{x^hy^k: 0\le h\le \beta-1, 0\le k\le \gamma-1\}.
$$
Then ([\PILD]), for $M=M(\beta,\gamma)$,
$$
D=\beta\gamma,\,\,\,\,\,\,\,
R={D\,(\gamma+\beta-2)\over 2},\,\,\,\,\,\,\,
S=D(\beta-1+\gamma-1)=2R,\,\,\,\,\,\,\,
C\le 2D,
$$
and (elementarily)
$$
\max\Big({1\over \beta},{1\over\gamma}\Big)\,\le\,\rho\,\le\,{1\over
\beta}+{1\over\gamma}.
$$

\medskip
\noindent
{\genb Proof of 1.4.\/} The proof adapts the proof of [\PILD, 1.4]
using 2.1 instead of [\PILD, 4.2].

\smallskip

Consider first a $C^\infty$ function $f$ on a subinterval of $[-1,1]$
with $|f'|\le 1$, with $f^{(j)}$ either nonvanishing or identically
vanishing for $j=0,\ldots,D$. Suppose $f$ satisfies an
irreducible algebraic relation of degree $(b,c), d=\max(b,c)$.
If $d=b$ take $M=(d,\delta)$ with $\delta\ge d$; if $d=c$
take $M=M(\delta,d)$ with $\delta\ge d$.
Then, by 2.1, $X(\QQ, H)$ is contained in the union of at most
$$
10 d^2\delta^2 4^d 2^{\rho} H^{2\rho}\le 20 d^2\delta^2 4^d 
H^{2/d+2/\delta}
$$
curves defined in $M$. The intersections are proper, $X$ is of degree
$\le b+c\le 2d$, the curves in $M$ of degree $\le 2\delta$, so
$$
N(X,H)\le 80 d^3\delta^3 4^d 
H^{2/d+2/\delta}.
$$

Next consider an algebraic curve $X$ defined by $F(x,y)=0$ in the box
$B=[-1,1]^2$, where $F$ is irreducible of bidegree $(b,c)$ and
$d=\max(b,c)$. Then $X$ has at most $2d(2d-1)$ singular points, and  at
most $4d(d-1)$ points with slope $\pm 1$. So $X\cap B$ consists of at
most $20d^3$ graphs of $C^\infty$ functions $f$ with slope $|f'|\le 1$
relative to one of the coordinate axes.

\smallskip

For each such function, the domain can be divided into at most
$8d^2D^2$ subintervals (see [\BOPI, Lemmas 5 and 6]) in which 
$f^{(j)}$ is nonvanishing or identically zero, $j=1,2,\ldots,D$. So
$$
N(X,H)\le 25.2^{10} d^{1o}\delta^5 4^d 
H^{2/d+2/\delta}.
$$

\smallskip

Finally, let $F(x,y)$ of bidegree $(b,c), d=\max(b,c)$,
$X=\{(x,y)\in\RR^2: F(x,y)=0\}$. Let
$P=(x,y)\in X(\QQ)$ with $H(P)\le H$. Then one of the following holds:

\smallbreak

(i) $|x|, |y| \le 1$

(ii) $|x|\le 1, |y|>1$

(iii) $|x|> 1, |y|\le 1$

(iv) $|x|> 1, |y| > 1$.

\smallskip

In case (i), $P$ lies in the box $[-1,1]^2\subset\RR^2$.
In case (ii), the point $Q=(x,1/y)$ is on the curve $Y: y^cF(x,1/y)=0$.
This curve is also irreducible and of bidegree $(b,c)$ (because $F$
must have a term independent of $y$). The point $Q$ is then in the box
$[-1,1]^2$ and has $H(Q)\le H$. Likewise in cases (iii) and (iv) the
corresponding points $R=(1/x,y), S=(1/x, 1/y)$ lie on irreducible
curves $x^bF(1/x,y)=0, x^b y^c F(1/x,1/y)=0$ of bidegree
$(b,c)$ in the box $[-1,1]^2$ and have height $\le H$.

\smallskip

Therefore, up to a factor 4, it suffices to consider the
points of $F$ inside the box $[-1,1]^2$, so
$$
N(X,H)\le 100 (2d)^{10} 4^d\delta^5 
H^{2/d+2/\delta}.
$$
Take $\delta$ to be the least integer exceeding $\log H$. Then, provided
$H\ge e^d$ (so that $\delta\ge d$), 
$$
N(X,H)\le 100 e^2 2^{15}\,d^{10}(\log H)^5 4^d 
H^{2/d}.
$$
However for $\log H\le d$ the bound is easily seen to hold as
well.\ \qed

\smallskip

\bigbreak
\centerline{\genb Acknowledgements}
\medskip

This paper was written while I was a visitor at the Mathematical
Institute, Oxford. I am grateful to the Institute, and in particular to
D. R. Heath-Brown and A. Lauder, for their hospitality. My stay in
Oxford was supported by my home institution, McGill University, and in
part by a grant from NSERC, Canada.

\vfil
\eject

\centerline{\genb References}

\bigskip

\item{\BOMB.} E. Bombieri, email 2003.

\item{\BOPI.} E. Bombieri and J. Pila, The number of 
integral points on arcs and ovals, {\genit Duke Math. J.\/}
{\genb 59\/} (1989), 337--357.

\item{\GAVO.} A. Gabrielov and N. Vorobjov, Complexity of computations
with pfaffian and noetherian functions, in {\genit Normal Forms,
Bifurcations and Finiteness problems in Differential Equations,\/}
Kluwer, 2004.

\item{\GWKP.} J. Gwozdziewicz, K. Kurdyka, and A. Parusinski,
On the number of solutions of an algebraic equation on the curve
$y=e^x+\sin x, x>0$, and a consequence for o-minimal structures,
{\genit Proc. Amer. Math. Soc.} {\genb 127} (1999), 1057--1064.

\item{\HEAT.} D. R. Heath-Brown, The density of rational
points on curves and surfaces, {\genit Ann. Math.\/} {\genb
155\/} (2002), 553--595.

\item{\KHOV.} A. G. Khovanskii, {\genit Fewnomials,\/} Translations of
Mathematical Monographs {\genb 88}, AMS, Providence, 1991.

\item{\LANG.} S. Lang, {\genit Introduction to transcendental
numbers,\/} Addison-Wesley, Reading, 1966.

\item{\PILA.} J. Pila, Geometric postulation of a smooth
function and the number of rational points, {\genit Duke
Math. J. }{\genb 63} (1991), 449--463.

\item{\PILC.} J. Pila, Integer points on the dilation of a subanalytic
surface, {\genit Quart. J. Math.\/} {\genb 55} (2004), 207--223

\item{\PILD.} J. Pila, Rational points on a subanalytic
surface, submitted.

\item{\SHID.} A. B. Shidlovskii, {\genit Transcendental numbers,}
translated from the Russian by N. Koblitz with a foreword by W. D.
Brownawell. de Gruyter Studies in Mathematics {\genb 12\/}, Walter de
Gruyter
\& Co., Berlin, 1989.

\bigskip
\bigskip
\bigskip

\line{Department of Mathematics and Statistics \hfil
Mathematical Institute}
\line{McGill University \hfil
University of Oxford}
\line{Burnside Hall\hfil
24-29 St Giles}
\line{805 Sherbrooke Street West\hfil
Oxford OX1 3LB}
\line{Montreal, Quebec, H3A 2K6\hfil
UK}
\leftline{Canada}

\bigskip
\bigskip

\leftline{pila@math.mcgill.ca}

\bigskip
\bigskip
\noindent
(Submitted for publication on 11 August 2004)

\bye